\documentclass[12pt]{article}

\usepackage[english]{babel}
\usepackage{amssymb,a4wide}
\usepackage{amsmath,xypic,epsf,mathrsfs,theorem}

\usepackage[all]{xy}
\newdir{ (}{{}*!/-5pt/\dir^{(}}
\newdir{ )}{{}*!/-5pt/\dir_{(}}

\newtheorem{theorem}{Theorem}[subsection]
\newtheorem{lemma}[theorem]{Lemma}
\newtheorem{corollary}[theorem]{Corollary}

\newtheorem{proposition}[theorem]{Proposition}

{\theorembodyfont{\upshape}

\newtheorem{remark}[theorem]{Remark}

}

\numberwithin{equation}{section}
\numberwithin{theorem}{section}

\newcommand{\Id}{\mathrm{id}}

\newcommand{\cO}{{\cal O}}

\newcommand{\cZ}{{\cal Z}}
\newcommand{\C}{{\mathbb C}}
\newcommand{\Z}{{\mathbb Z}}

\newcommand{\Cs}{{$C^*$-al\-ge\-bra}}

\newcommand{\sh}{{$^*$-ho\-mo\-mor\-phism}}

\newenvironment{proof}[1][Proof:]
{\begin{trivlist}\item[]\textbf{#1} }
{\hbox{}\nobreak\hfill\quad\hbox{$\square$}\end{trivlist}}

\begin{document}
\title{Strongly Self-Absorbing \Cs s
  which contain a nontrivial projection}

\author{Marius Dadarlat\footnote{M.D.\ was partially supported by NSF grant
  \#DMS-0801173.} \, and Mikael R\o rdam\footnote{ M.R.\ was supported by
  a grant from the Danish Natural Science Research Council (FNU)}}
\date{}

\maketitle

\begin{center} \it{Dedicated to Joachim Cuntz on the occasion of his
    60th birthday}

\end{center}

\begin{abstract} \noindent
It is shown that a strongly self-absorbing \Cs{} is
  of real rank zero and absorbs the Jiang-Su algebra if it contains a
  non-trivial projection. We also consider cases where the UCT is
  automatic for strongly self-absorbing \Cs s, and $K$-theoretical ways of
  characterizing when Kirchberg algebras are strongly self-absorbing.
\end{abstract}

\section{Introduction} 

\noindent Strongly self-absorbing \Cs s were first systematically
studied by Toms and Winter in \cite{TW}. The
classfication program of Elliott had prior to that been seen to work
out particulaly  
well for those (separable, nuclear) \Cs s that tensorially absorb one
of the Cuntz algebras $\cO_2$, $\cO_\infty$, or the Jiang-Su algebra
$\cZ$. More precisely, thanks to deep theorems of Kirchberg, 
the classification of separable, nuclear, stable \Cs
s that absorb the Cuntz algebra $\cO_2$ is complete (the invariant is
the primitive ideal space); and separable, nuclear, stable \Cs s
that absorb the Cuntz algebra $\cO_\infty$ are classified by an ideal
related $KK$-theory. The situation for separable, nuclear \Cs s that
absorb the Jiang-Su algebra  is at present
very promising (see for example \cite{W}) but not as complete as in
the purely infinite case.

The \Cs s $\cO_2$, $\cO_\infty$ and $\cZ$ are all examples of strongly
self-absorbing \Cs s. They are in \cite{TW} defined to be those unital
separable \Cs s $D\neq \mathbb{C}$ for which there is an isomorphism
$D \to D \otimes D$ that 
is approximately unitarily equivalent to the \sh{} $d \mapsto d
\otimes 1$. Strongly self-absorbing \Cs s are automatically simple and
nuclear, and they have at most one tracial state. 
It is shown in \cite{TW} that if $D$ is a strongly self-absorbing
\Cs{} in the UCT class, then it has the same $K$-theory as one of the
\Cs s in the following list: $\cZ$, UHF-algebras of infinite type,
$\cO_\infty$, $\cO_\infty$ tensor a UHF-algebra of infinite type, or
$\cO_2$. It is an open problem if nuclear \Cs s always satisfy the UCT
(and also if strongly self-absorbing \Cs s enjoy this property);
and it is an intriguing problem, very much related to the Elliott
classification program, if the list above exhausts all
strongly self-absorbing \Cs s.  Should the latter
be the case, 
then it would in particular follow that every strongly self-absorbing
\Cs{} absorbs the Jiang-Su algebra $\cZ$. 
By the Kirchberg-Phillips
classification theorem, a strongly self-absorbing Kirchberg
algebra belongs to the list above if and only if it belongs to the
UCT class. Let us also remind the reader that a strongly
self-absorbing \Cs{} is a Kirchberg algebra if and only if it is not
stably finite (or, equivalently, if and only if it is traceless).

In Section~2 of this paper we show that every strongly self-absorbing
\Cs{} which 
contains a non-trivial projection is of real rank zero and absorbs the
Jiang-Su algebra. In Section~3 we consider $K$-theoretical
conditions on strongly self-absorbing Kirchberg algebras. One such
condition (phrased at the level of $K$-homology) characterizes the
Kirchberg algebra $\cO_\infty$, and other results in Section~3 give
$K$-theoretical characterizations on when a Kirchberg algebra is
strongly self-absorbing. 

\section{Strongly self-absorbing \Cs s with a non-trivial projection}

In this section we show that any strongly self-absorbing \Cs{} that
contains a non-trivial projection is automatically approximately
divisible, of real rank zero, and absorbs the Jiang-Su algebra $\cZ$

\begin{lemma} \label{lm:0} There is a unital \sh{} from
  $M_3 \oplus M_2$ into a unital \Cs{} $A$ if and only if $A$ contains
  projections 
  $e, e'$ such that $e \perp e'$, $e \sim e'$, and $1-e-e' \precsim e$.
\end{lemma}

\begin{proof} It is easy to see that such projections $e$ and $e'$
  exist in $M_3 \oplus M_2$ and hence in any unital \Cs{} $A$ that
is the target of a unital \sh{} from $M_3 \oplus M_2$.

Assume now that such projections $e$ and $e'$ exist. 
Let $v \in A$ be a partial isometry such that $v^*v = e$ and $vv^* =
e'$. Put $f_0 = 1-e-e'$. Find a subprojection $f_1$ of $e$ which is
equivalent to $f_0$, and put $f_2 = vf_1v^*$. Put $g_1 = e-f_1$ and
put $g_2 = e'-f_2 = vg_1v^*$. 
The projections $f_0,f_1,f_2,g_1,g_2$ then satisfy
$$1 = f_0+f_1+f_2+g_1+g_2, \qquad f_0 \sim f_1 \sim f_2, \qquad g_1
\sim g_2.$$
Extending the sets $\{f_0,f_1,f_2\}$ and $\{g_1,g_2\}$ to sets of
matrix units for $M_3$ and $M_2$, respectively, yields a unital \sh{}
from $M_3 \oplus M_2$ into $A$. (If the $f_j$'s are zero or if the
$g_j$'s are zero, then this \sh{} will fail to be injective, and will
instead give a unital embedding of $M_2$ or $M_3$ into $A$.)
\end{proof}

\noindent If $D$ is any unital \Cs{} then we let $D^{\otimes n}$ denote the
$n$-fold tensor product $D \otimes D \otimes \cdots \otimes D$ (with
$n$ tensor factors), and we let $D^{\otimes \infty}$ denote the
infinite tensor product $\bigotimes_{n=1}^\infty D$. The latter is the
inductive limit of the sequence
$$D \to D^{\otimes 2} \to D^{\otimes 3} \to D^{\otimes 4} \to
\cdots,$$
(with connecting mappings $d \mapsto d \otimes 1_D$).
We shall view $D$ as a (unital) sub-\Cs{} of $D^{\otimes n}$,
$D^{\otimes n}$ as a sub-\Cs{} of $D^{\otimes m}$ (if $n \le m$), and
finally $D$ and $D^{\otimes n}$ are viewed as subalgebras of
$D^{\otimes \infty}$.

If $x \in D^{\otimes n}$, then $x^{\otimes k}$ will denote the
$k$-fold tensor product
$$x^{\otimes k} = x \otimes x \otimes x \otimes \cdots \otimes x \in
D^{\otimes kn}.$$

The proof of the lemma below resembles the proof of 
\cite[Lemma~6.4]{RW}. 

\begin{lemma} \label{lm1} Let $D$ be a strongly self-absorbing \Cs,
  and let $p$ be a projection in $D$. Consider the
  following projections in $D \otimes D$,
$$e_1 = p \otimes (1-p), \qquad e_1' = (1-p) \otimes p, \qquad f = p
\otimes p + (1-p) \otimes (1-p).$$
For each natural number $n$ consider also the following
projections in $D^{\otimes 2(n+1)}$,
$$e_{n+1} = f^{\otimes n} \otimes p \otimes (1-p), \qquad
 e_{n+1}' = f^{\otimes n} \otimes (1-p) \otimes p.$$
It follows that the projections $e_1,e_2, \dots, e_1',e_2', \dots$ are
pairwise orthogonal in $D^{\otimes \infty}$,
and that $e_j \sim e_j'$. Moreover, for each
natural number $n$, set
$$E_n = e_1+e_2+ \cdots + e_n, \qquad E_n' = e_1' + e_2' + \cdots +
e_n'.$$
Then $E_n \perp E_n'$, $E_n \sim E_n'$, and
\begin{equation} \label{eq:1}
1 - E_n - E_n' = f^{\otimes n}.
\end{equation}
\end{lemma}

\begin{proof} The equivalence $e_j \sim e_j'$ comes from the fact that
  the flip automorphism $a \otimes b \mapsto b \otimes a$ on $D
  \otimes D$ is approximately inner when $D$ is strongly
  self-absorbing. The projections $e_1,e_2, \dots,
  e_1',e_2', \dots$ are pairwise orthogonal by construction. 
  The only thing left to
  prove is \eqref{eq:1}. We prove this by induction after $n$, and
  note first that \eqref{eq:1} for $n=1$ follows from the fact that
  $e_1+e'_1+f=1$. Suppose that 
  \eqref{eq:1} holds for some $n \ge 1$. Then 
\begin{eqnarray*}
1-E_{n+1}-E_{n+1}' & = & 1 - E_n - E_n' - e_{n+1}-e_{n+1}' \\
&=& f^{\otimes n} \otimes 1_D \otimes 1_D - f^{\otimes n} \otimes p
\otimes (1-p) - f^{\otimes n} \otimes (1-p) \otimes p \\
&=& f^{\otimes (n+1)}.
\end{eqnarray*}
\end{proof}

\begin{lemma} \label{lm:2}
Let $D$ be a strongly self-absorbing \Cs{} and let $p$ be a
projection in $D$ such that $p \ne 1$. Then there exists a natural number $n$
such that $p^{\otimes n} \precsim 1 - p^{\otimes n}$ in $D^{\otimes n}$.
\end{lemma}

\begin{proof} To simplify the notation we express our calculations in
  terms of the monoid $V(D)$ of Murray-von Neumann equivalence classes
  of projections in $D$ and in matrix algebras over $D$. Let $[e] \in
  V(D)$ denote the equivalence class containing the projection $e$ in
  (a matrix algebra over) $D$.

Since $D$ is simple and $p \ne 1$ 
there is a natural number $n$ such that $n[1-p] \ge
[p]$. It follows that
\begin{eqnarray*}
[1 - p^{\otimes n}] & \ge & [(1-p) \otimes p \otimes \cdots \otimes p]
+ [p \otimes (1-p) \otimes \cdots \otimes p] + [p \otimes p \otimes
\cdots \otimes (1-p)] \\
&=& n[(1-p) \otimes p \otimes \cdots \otimes p] \\
&\ge & [p \otimes p \otimes p \otimes \cdots \otimes p]
= [p^{\otimes n}],
\end{eqnarray*}
where the equality between the second and third expression holds
because the flip on a strongly self-absorbing \Cs{} is approximately
inner. 
\end{proof}

\begin{lemma} \label{lm:3}
Let $D$ be a strongly self-absorbing \Cs{}, let $p$ be a
  projection in $D^{\otimes k}$, and let $e$ be a projection in
  $D^{\otimes l}$ for some natural numbers $k$ and $l$. Assume that $p
  \ne 1$ and that $e \ne 0$. It follows that there exists a natural
  number $n$ such that $p^{\otimes n} \precsim e$ in $D^{\otimes
    \infty}$. 
\end{lemma}

\begin{proof} 
Let $d$ be a natural number such that $dk \ge l$. Upon replacing
$p$ with $p^{\otimes d}$, $e$ with $e \otimes 1_D^{\otimes (dk-l)}$, and $D$
with $D^{\otimes dk}$ we can assume that $p$ and $e$ both belong to
$D$. Use Lemma~\ref{lm:2} to find $m$ such that $p^{\otimes m} \precsim
1-p^{\otimes m}$. By replacing $p$ with $p^{\otimes m}$, $e$ with
$e \otimes 1_D^{\otimes (m-1)}$, and $D$ 
with $D^{\otimes m}$ we can assume that $p$ and $e$ both belong to
$D$ and that $p \precsim 1-p$. 

Now, $p \sim q \le 1-p$ for some
projection $q$ in $D$. In the language of the monoid
$V(D)$ we have
$$[{1_D}^{\otimes k}] \ge [(p + q)^{\otimes k}] = 2^k[p^{\otimes k}]$$
for any natural number $k$.
Using simplicity of $D$ we can choose $n$ such that $2^{n-1}[e] \ge [p]$. Then
$$[e] = [e \otimes {1_D}^{\otimes (n-1)}] \ge 2^{n-1}[e \otimes p^{\otimes
  (n-1)}] \ge [p^{\otimes n}],$$
in $V(D^{\otimes n})$ as desired, where we in the first identity
have used that the embedding of $D$ into $D^{\otimes n}$ maps $e$ onto
$e \otimes {1_D}^{\otimes (n-1)}$. 
\end{proof}

\begin{theorem} \label{thm:1}
Let $D$ be a strongly self-absorbing \Cs. Then the following three
conditions are equivalent:
\begin{enumerate}
\item $D$ contains a non-trivial projection (i.e., a projection other
  than $0$ and $1$).
\item $D$ is approximately divisible.
\item $D$ is of real rank zero.
\end{enumerate}
If any of the three equivalent conditions are satisfied, then $D$
absorbs the Jiang-Su algebra, i.e., $D \cong D \otimes \cZ$.
\end{theorem}

\begin{proof} (i) $\Rightarrow$ (ii). If $D$ is strongly
  self-absorbing, then there is an asymptotically central sequence of
  embeddings of $D$ into itself, i.e., a sequence $\rho_k \colon D \to
  D$ of unital \sh s such that $\|\rho_k(x)y - y\rho_k(x)\| \to 0$ as
  $k\to \infty$ for all $x,y \in D$. 

Identify $D$ with $D_0^{\otimes \infty}$ where $D_0 \cong D$. Take a
non-trivial projection $p$ in $D_0$. For each
natural number $n$ let $E_n, E_n' \in D_0^{\otimes 2n}$ be as in
Lemma~\ref{lm1} (corresponding to our non-trivial projection $p$). 
Then $e_n\neq 0, E_n\neq 0, $ and so $0\neq f^{\otimes n}\neq 1$.
Use \eqref{eq:1} and Lemma~\ref{lm:3} to find $n$
such that $1-E_n - E_n' \precsim p \otimes (1-p) \le E_n$. It
then follows from Lemma~\ref{lm:0} that there is an injective unital \sh{} from
$M_3 \otimes M_2$ into $D_0^{\otimes 2n} \subseteq D$. Composing this
unital \sh{} with the unital \sh{} $\rho_k$ yields an asymptotically
central sequence of unital \sh s from $M_3 \otimes M_2$ into $D$. This
shows that $D$ is approximately divisible.

(ii) $\Rightarrow$ (iii). It is shown in \cite{BKR:ApprDiv} that a
simple approximately divisible \Cs{} is of real rank zero if and only
if projections in the \Cs{} separate the quasitraces. As quasitraces
on a exact \Cs{} are traces, \cite{Haa:quasitraces}, a result that
applies to our case since strongly self-absorbing \Cs s are nuclear
and hence exact, and since a strongly self-absorbing \Cs{} has at most
one tracial state, quasitraces are automatically separated by just one
projection, say the unit. 

(iii) $\Rightarrow$ (i). This is trivial. The only \Cs{} of real rank
zero that does not have a non-trivial projection is $\C$, the algebra
of complex numbers. This \Cs{} is not strongly self-absorbing by
convention.

Finally, any simple approximately divisible \Cs{} is $\cZ$-absorbing,
cf.\ \cite{TW}.
\end{proof}

\begin{lemma}
 Let $D$ be a strongly self-absorbing \Cs.  Then $K_0(D)$ has a
 natural structure 
of commutative unital ring with unit $[1_D]$. If $\tau$ is a unital trace on $D$, then $\tau$ induces a morphism
of unital rings $\tau_* \colon K_0(D)\to \mathbb{R}$.
\end{lemma}
\begin{proof}
 Fix an isomorphism $\gamma \colon
D\otimes D \to D$. The multiplication on $K_0(D)$
is defined by composing $\gamma_* \colon
K_0(D\otimes D) \to K_0(D)$ with the canonical map 
 $K_0(D)\otimes K_0(D) \to K_0(D\otimes D)$. Since any two unital
\sh{}
 from $D\otimes D$ to $D$ are approximately unitarily equivalent, the
 above multiplication is well-defined and commutative. We leave the
 rest of proof for the reader, but note that if $D$ has a unital
 trace, then $\tau \otimes \tau$ is 
the unique unital trace of $D\otimes D$.  
\end{proof}
\begin{proposition} \label{prop:qd}
Let $D$ be a strongly self-absorbing \Cs. 
Suppose  that $D$ is quasidiagonal and that $K_0(D)$ is torsion free.
Then either $K_0(D)\cong \mathbb{Z}$ or there is a \emph{UHF} algebra
$B$ of infinite type such that $K_0(D)\cong K_0(B)$. If, in addition,
$D$ is assumed to contain a nontrivial projection, then $D\otimes B
\cong D$, where $B$ is as above. 
\end{proposition}
\begin{proof} Since $D$ is quasidiagonal it embeds unitally in the
  universal UHF algebra $B_{\mathbb{Q}}$ and $D\otimes B_{\mathbb{Q}}
  \cong  B_{\mathbb{Q}}$, as explained in \cite[Rem.\ 3.10]{DW}. 
The restriction of the unital trace of $B_{\mathbb{Q}}$ to $D$ is
denoted by $\tau$. 
Thus we have an exact sequence
\begin{equation*}
    \xymatrix{
0\ar[r] & H \ar[r] & K_0(D) \ar[r]^-{\tau_*} & \tau_* K_0(D) \ar[r] & 0
}
\end{equation*}
where $H$ is the kernel of $\tau_*$.  Since $\mathbb{Z}\subseteq \tau_*
K_0(D) \subseteq \mathbb{Q}$, 
and $K_0(D)\otimes \mathbb{Q}\cong\mathbb{Q}$, the map 
$\tau_*\otimes \mathrm{id}_{\mathbb{Q}} \colon 
K_0(D) \otimes \mathbb{Q} \to \tau_* K_0(D) \otimes \mathbb{Q}$ is an
isomorphism. Therefore $H\otimes \mathbb{Q}=0$ and so $H$ is a torsion
subgroup of $K_0(D)$. But we assumed that $K_0(D)$ is torsion free and
hence $H=\{0\}$ and 
$\tau_* \colon K_0(D)\to \tau_* K_0(D)\subseteq \mathbb{Q}$ is an
isomorphism of unital rings. 
The unital subrings of $\mathbb{Q}$ are easily determined and
well-known. They are 
parametrized by arbitary sets $P$ of prime numbers. For each $P$ the
corresponding ring  $R_P$ 
 consists of rational numbers $r/s$ with $r$ and $s$ relatively prime
 and such that all prime factors of $s$ are in $P$. If $P=\emptyset$
 then $R_P=\mathbb{Z}$,  otherwise $R_P$ is 
isomorphic to the $K_0$-ring associated to a UHF algebra $B$ of infinite type.

Suppose now that $D$ contains a nontrivial projection. By 
Theorem~\ref{thm:1}, $D$ has real rank zero and absorbs the Jiang-Su
algebra 
$\cZ$. In particular, $K_0(D)$
is not $\Z$ and is hence isomorphic (as a scaled abelian group) to
$K_0(B)$ for some UHF-algebra $B$ of 
infinite type. It follows from \cite{Ror:Z} that $D$ has 
stable rank one and that $K_0(D)$ is weakly unperforated. Moreover, by
\cite[Sect.\ 6.9]{Blackadar'sbook}, $K_0(D)$ has the strict order
induced by $\tau_*$. The isomorphism $K_0(B)\cong K_0(D)$ of
scaled abelian groups is therefore an order isomorphism, and by 
the properties of $D$ established above we can conclude that $B$
embeds unitally into $D$, whence $D\otimes B \cong D$.
\end{proof}


\begin{corollary} 
Let $D$ be a strongly self-absorbing \Cs\ with torsion free $K_0$-group.
Suppose that $D$ contains a non-trivial projection and that $D$ embeds unitally
into the \emph{UHF} algebra $M_{p^\infty}$ for some prime number $p$.
Then $D\cong M_{p^\infty}$.
\end{corollary}
\begin{proof}
 By Proposition~\ref{prop:qd} there is a prime $q$ such that
 $M_{q^\infty}$ in contained unitally in $D$ and hence in
 $M_{p^\infty}$. From this we deduce that $q=p$. Finally since 
$M_{p^\infty}\subseteq D \subseteq M_{p^\infty}$ we conclude that $D\cong
M_{p^\infty}$. 
\end{proof}

\section{Strongly self-absorbing algebras and K-theory}
The class of strongly self-absorbing Kirchberg algebras satisfying the
UCT was completely described in \cite{TW}. In this section we give
properties and characterizations of strongly self-absorbing Kirchberg
algebras which can be derived without assuming the UCT. 
For unital a \Cs\ $D$ we denote by $\nu_D$ the unital \sh{}
$\mathbb{C}\to D$. 
When the \Cs\ $D$ is clear from context we will write $\nu$ instead of $\nu_D$.

\begin{proposition}\label{oinf} Let $D$ be a  strongly self-absorbing 
\Cs. If $D$ is not finite and the unital \sh{} $\mathbb{C}\to D$
induces a surjection $K^0(D)\to K^0(\mathbb{C})$, then $D\cong O_\infty$.
\end{proposition}

\begin{proof} 
By \cite[Prop.\ 5.12]{TW}, two  strongly self-absorbing
\Cs s are isomorphic if and only if they embed unitally into each other.
Thus it suffices to show the existence of unital \sh s
$\cO_\infty\to D$  and $D\to \cO_\infty$ 
Since $D$ is not finite, it must be a Kirchberg algebra, see
\cite[Sec.\ 1]{TW},  
and hence
$\cO_\infty$ embeds unitally in $D$ by \cite[Prop.\ 4.2.3]{Ror:ency}. 
It remains to show that $D$ embeds unitally
in $\cO_\infty$. 

By assumption, the map $\nu^* \colon
KK(D,\mathbb{C})\to KK(\mathbb{C},\mathbb{C})$
is surjective. By multiplying with the $KK$-equivalence class
given by the unital morphism $\mathbb{C}\to \cO_\infty$, we obtain
that the map
$\nu^* \colon KK(D,\cO_\infty)\to KK(\mathbb{C},\cO_\infty)$
is surjective. If $\varphi \colon D \to \cO_\infty \otimes \mathcal{K}$
is a \sh, then, after identifying
$KK(\mathbb{C},\cO_\infty)\cong K_0(\cO_\infty)$, the map $\nu^*$ sends
$[\varphi]$ to the class 
$[\varphi(1_D)]\in K_0(\cO_\infty)$. By \cite[Thm.\ 8.3.3]{Ror:ency}
each element of $KK(D,\cO_\infty)$ is represented by a
$^*$-homomorphism. 
Therefore, by the surjectivity
of $\nu^*$, there is a $^*$-homomorphism $\varphi \colon 
D \to \cO_\infty \otimes \mathcal{K}$ such that 
$[\varphi(1_D)]=[1_{\cO_\infty}]$. Since these are both  full
projections, by \cite[Prop.\ 4.1.4]{Ror:ency} there is a partial
isometry $v\in \cO_\infty \otimes \mathcal{K}$ such that 
$v^*v= \varphi(1_D)$ and $v v^*=1_{\cO_\infty}$. Then $v \varphi v^*$
is a unital 
embedding $D\to \cO_\infty$.
\end{proof}
\begin{remark}
 Note that the isomorphism $D\cong \cO_\infty$
 was obtained without 
assuming that $D$ satisfies the UCT. Let us argue that assumptions of
Proposition~\ref{oinf} are natural.  
Let $A$ and $B$ be unital \Cs s and let $\nu \colon
\mathbb{C}\to A$ and $\nu \colon \mathbb{C}\to B$
be the corresponding unital $^*$-homomorphisms.
The condition that there is a morphism of pointed groups
$(K_0(A),[1_A])\to (K_0(B),[1_B])$ 
can be viewed as the condition that the diagram

\begin{equation*}
    \xymatrix{
A & & B\\
& \mathbb{C}\ar[ur]_\nu\ar[ul]^\nu
}
\end{equation*}
can be completed to a commutative diagram after passing to $K$-theory:

\begin{equation*}\label{unit1}
    \xymatrix{
K_0(A)\ar@{-->}[rr] & & K_0(B)\\
& K_0(\mathbb{C})\ar[ur]_{\nu_*}\ar[ul]^{\nu_*}
}
\end{equation*}

It would then be completely natural to use $K$-homology instead of $K$-theory
and ask that the first diagram  can be completed to a commutative
diagram after passing 
to $K$-homology.

\begin{equation*}\label{unit2}
    \xymatrix{
K^0(A)\ar[dr]_{\nu^*} & &\ar@{-->}[ll] K^0(B)\ar[dl]^{\nu^*}\\
& K^0(\mathbb{C})
}
\end{equation*}

Now let us observe that the condition, imposed in
Proposition~\ref{oinf}, that $\nu^* \colon 
K^0(D)\to  K^0(\mathbb{C})$ is surjective  clearly is equivalent to the
existence of a commutative  diagram 
\begin{equation*}\label{unit3}
    \xymatrix{
K^0(\cO_\infty)\ar[dr]_{\nu^*} & &\ar[ll]_-{\alpha} K^0(D)\ar[dl]^{\nu^*}\\
& K^0(\mathbb{C})
}
\end{equation*}
where $\alpha$ is a surjective morphism.

If $D$ satisfies the UCT, then the condition above can be translated 
in terms of K-theory as follows. Since the commutative diagram
\begin{equation*}
    \xymatrix{
K^0(D)\ar[d]_{\nu^*}\ar[r]& \mathrm{Hom}(K_0(D), \mathbb{Z})\ar[d]\\
K^0(\mathbb{C})\ar[r]& \mathrm{Hom}(K_0(\mathbb{C}), \mathbb{Z})
}
\end{equation*}
has surjective horizontal arrows,
 the assumption on $K$-homology in Proposition~\ref{oinf}
is equivalent  the existence a  group homomorphism $K_0(D)\to
\mathbb{Z}$ which maps $[1_D]$ to $1$. This is obviously equivalent to
the condition that $[1_D]$ is an infinite order element of $K_0(D)$
and that the subgroup 
that it generates, $\mathbb{Z}[1_D]$, is a direct summand of $K_0(D)$.
\end{remark}

Our next goal is to show that for a unital Kirchberg algebra the
property of being strongly self-absorbing is purely a $KK$-theoretical
condition. 
Let 
$$C_\nu=\{f \colon [0,1]\to D \mid f(0)\in \mathbb{C}1_D, \quad f(1)=0\}$$
be the mapping cone of the unital \sh s $\nu \colon \mathbb{C}\to D$.
\begin{proposition}\label{map_cone} Let $D$ be a unital  Kirchberg
  algebra. Then $D$ is strongly self-absorbing if and only if $C_{\nu}
  \otimes D$ is $KK$-equivalent to zero. 
\end{proposition}
\begin{proof} 
We begin with a general observation.
For a  \sh{}  $\varphi \colon A \to B$ of separable 
  \Cs s and any separable \Cs\ $C$, there is an exact Puppe
 sequence in  $KK$-theory  
(\cite[Thm.\ 19.4.3]{Blackadar'sbook}): 
\begin{equation*}
    \xymatrix{
KK(B, C) \ar[r]^{\varphi^*} &KK(A,C) \ar[r] & KK(C_\varphi,C)\ar[d]\\
KK_1(C_\varphi,C )\ar[u] &\ar[l] KK_1(A,C) & \ar[l]^{\varphi^*} KK_1(B,C) 
}
\end{equation*}
It is apparent that $[\varphi]\in
KK(A,B)^{-1}$ if and only if composition with $[\varphi] \in KK(A,B)$
induces a bijection 
$\varphi^* \colon KK(B,C)\to KK(A,C)$  for any separable \Cs{} $C$, 
or equivalently, for just  $C=A$ and $C=B$. Therefore, by the
exactness of the Puppe sequence, we see that 
that $\varphi$ induces a  $KK$-equivalence if and only if
its mapping cone \Cs{} $C_\varphi$ is $KK$-contractible. 

By applying this observation to the unital \sh{}
$\nu\otimes \Id_D \colon D \to D\otimes D$ 
we deduce that $\nu\otimes \Id_D$ induces a $KK$-equivalence if and only
if its mapping cone 
$C_{\nu \otimes \Id_D}\cong C_{\nu}\otimes D$ is 
$KK$-contractible.
Suppose now 
that $D$ is a strongly self-absorbing  Kirchberg algebra. 
Then $\nu\otimes \Id_D$ is asymptotically unitarily equivalent to a an
isomorphism by 
\cite[Thm.\ 2.2]{DW} and hence $\nu\otimes \Id_D$ induces a
$KK$-equivalence. Conversely,  
if $\nu\otimes \Id_D$ induces a $KK$-equivalence, then $\nu\otimes \Id_D$
is  asymptotically unitarily equivalent to an isomorphism $D\to
D\otimes D$ by \cite[Thm.\ 8.3.3]{Ror:ency} 
and hence $D$ is strongly self-absorbing.
\end{proof}

\noindent
We have the following result related to Proposition~\ref{map_cone}.

\begin{proposition}\label{mappping_cone} Let $D$ be a unital Kirchberg
  algebra such that $D\cong D\otimes D$.  
The following assertions are equivalent:
\begin{enumerate}
\item $D$  is strongly self-absorbing.
\item $KK(C_\nu,SD)=0$.
\item $KK(C_\nu,D\otimes A)=0$ for all separable \Cs s $A$.
\item The map $KK(D,D\otimes A)\to KK(\mathbb{C},D\otimes A)$
 is bijective for  all separable \Cs s $A$.
\end{enumerate}
\end{proposition}

\begin{proof} (iii) $\Leftrightarrow$ (iv). This equivalence  is
  verified by using the Puppe sequence associated to $\nu \colon 
\mathbb{C} \to D$, arguing as in the proof of Proposition~\ref{map_cone}.

(i) $\Rightarrow$ (iv). 
This implication is proved in \cite[Thm.\ 3.4]{DW}.

(iii) $\Rightarrow$ (ii). This follows by taking $A=S\mathbb{C}$ in (iii).

(ii) $\Rightarrow$ (i). Fix an isomorphism $\gamma \colon D \to D\otimes D$.
Since $KK_1(C_\nu,D\otimes D)=0$ by hypothesis, it follows from the
Puppe sequence 
that the map 
$\nu^* \colon  KK(D,D\otimes D)\to KK(\mathbb{C},D\otimes D)$ is injective. 
 Therefore $\gamma$ and $\nu \otimes \Id_D$ induce the same class in
 $KK(D,D\otimes D)$ 
since they are both unital. It follows that $\nu \otimes \Id_D$ is
asymptotically unitarily 
equivalent to $\gamma$ and so $D$ is strongly self-absorbing.
\end{proof}

\begin{corollary} Let $D$ be a unital Kirchberg algebra such that
  $D\cong D\otimes D$. 
Then $D$ is strongly self-absorbing if and only if $\pi_2 \,
\mathrm{Aut}(D)=0$. 
\end{corollary}
\begin{proof}
 Since $\pi_2 \, \mathrm{Aut}(D)\cong KK(C_\nu,SD)$ by \cite[Cor.\ 3.1]{DW},
the conclusion follows from Proposition~\ref{mappping_cone}.
\end{proof}

\noindent
It was shown in \cite[Prop.\ 4.1]{Dad;K-theory} that if a unital
Kirchberg algebra satisfies the UCT, 
then $D$ is strongly self-absorbing if and only if the homotopy
classes $[X,\mathrm{Aut}(D)]$ reduces to singleton for any path
connected 
compact metrizable space $X$.

\vspace{.5cm}

\noindent{\sc Department of Mathematics, University of Copenhagen,
Universitetsparken 5, 2100 Copenhagen \O, Denmark}

\vspace{.3cm}

\noindent{\sl E-mail address:} {\tt rordam@math.ku.dk}\\
\noindent{\sl Internet home page:}
{\tt www.math.ku.dk/$\,\widetilde{\;}$rordam} \\

\vspace{.5cm}

\noindent{\sc Department of Mathematics, Purdue University 
West Lafayette, IN 47906, USA}

\vspace{.3cm}

\noindent{\sl E-mail address:} {\tt mdd@math.purdue.edu}\\
\noindent{\sl Internet home page:}
{\tt www.math.purdue.edu/$\,\widetilde{\;}$mdd} \\

\end{document}